\documentclass[12pt]{article}

\usepackage{latexsym}
\usepackage{amsfonts}
\usepackage{amssymb}
\usepackage{amsmath}
\usepackage{mathrsfs}
\input xypic
\usepackage{hyperref}

\newcommand{\sgn}{\mbox{sign\,}}

\newcommand{\be}{\begin{equation}}
\newcommand{\ee}{\end{equation}}
\newcommand{\bea}{\begin{eqnarray}}
\newcommand{\eea}{\end{eqnarray}}
\newcommand{\bean}{\begin{eqnarray*}}
\newcommand{\eean}{\end{eqnarray*}}
\newcommand{\brray}{\begin{array}}
\newcommand{\erray}{\end{array}}

\newtheorem{dfn}{Definition}[section]
\newtheorem{thm}[dfn]{Theorem}
\newtheorem{lmma}[dfn]{Lemma}
\newtheorem{ppsn}[dfn]{Proposition}
\newtheorem{crlre}[dfn]{Corollary}
\newtheorem{xmpl}[dfn]{Example}
\newtheorem{rmrk}[dfn]{Remark}

\newcommand{\bdfn}{\begin{dfn}\rm}
\newcommand{\bthm}{\begin{thm}}
\newcommand{\blmma}{\begin{lmma}}
\newcommand{\bppsn}{\begin{ppsn}}
\newcommand{\bcrlre}{\begin{crlre}}
\newcommand{\bxmpl}{\begin{xmpl}}
\newcommand{\brmrk}{\begin{rmrk}\rm}

\newcommand{\edfn}{\end{dfn}}
\newcommand{\ethm}{\end{thm}}
\newcommand{\elmma}{\end{lmma}}
\newcommand{\eppsn}{\end{ppsn}}
\newcommand{\ecrlre}{\end{crlre}}
\newcommand{\exmpl}{\end{xmpl}}
\newcommand{\ermrk}{\end{rmrk}}

\author{S. Sundar\\
Indian Statistical Institute, Delhi, \\ 
NewDelhi-110016.\\
\texttt{sundarsobers@gmail.com}}
\title{A computation with the Connes-Thom isomorphism}
\date{\today}

\begin{document}
\maketitle
\begin{abstract}
 Let $A \in M_{n}(\mathbb{R})$ be an invertible matrix. Consider the semi-direct product $\mathbb{R}^{n} \rtimes \mathbb{Z}$ where $\mathbb{Z}$ acts on $\mathbb{R}^{n}$ by matrix multiplication. Consider a strongly continuous action $(\alpha,\tau)$ of $\mathbb{R}^{n} \rtimes \mathbb{Z}$ on a $C^{*}$-algebra $B$ where $\alpha$ is a strongly continuous action of $\mathbb{R}^{n}$ and $\tau$ is an automorphism. The map $\tau$ induces a map $\widetilde{\tau}$ on $B \rtimes_{\alpha} \mathbb{R}^{n}$. We show that, at the $K$-theory level, $\tau$ commutes with the Connes-Thom map if $\det(A)>0$ and anticommutes if $\det(A)<0$. As an application, we recompute the $K$-groups of the Cuntz-Li algebra associated to an integer dilation matrix.
\end{abstract}
{\bf AMS Subject Classification No.:} {\large
46}L{\large 80}, {\large 58}B{\large 34}.  

\section{Introduction}
In \cite{Cuntz-Li}, Cuntz and Li initiated the study of  $C^{*}$-algebras associated with rings. In \cite{Cuntz}, Cuntz had earlier studied the $C^{*}$-algebra associated to the $ax+b$ group over $\mathbb{N}$. These $C^{*}$-algebras are unital, purely infinite and simple. Thus they are classified by their $K$-groups. In a series of papers, \cite{Cuntz-Li-1}, \cite{Cuntz-Li2} and \cite{wolf-Li}, the $K$-groups of these algebras were computed for number fields and for function fields. The main tool used in the $K$-group computation is the duality result proved in \cite{Cuntz-Li-1} and its variations.

Other approaches  and possible generalisations were considered in \cite{Ex2}, \cite{Quigg-Landstad} and in \cite{Sundar-Cuntz-Li}. The $C^{*}$-algebras studied in  \cite{Quigg-Landstad} and in \cite{Sundar-Cuntz-Li} were called Cuntz-Li algebras. Following \cite{Quigg-Landstad}, in \cite{Sundar-Cuntz-Li}, the Cuntz-Li algebra associated to a pair $(N \rtimes H,M)$ satisfying certain conditions was studied. Here $M$ is a normal subgroup of $N$ and $N \rtimes H$ is a semidirect product . The main example considered in \cite{Sundar-Cuntz-Li} is the Cuntz-Li algebra, denoted $\mathcal{U}_{\Gamma}$ associated to the pair $(\mathbb{Q}^{n} \rtimes \Gamma, \mathbb{Z}^{n})$ where $\Gamma$ is a subgroup of $GL_{n}(\mathbb{Q})$ acting by matrix multiplication on $\mathbb{Q}^{n}$. In \cite{Sundar-Cuntz-Li}, it was proved that $\mathcal{U}_{\Gamma}$ is Morita-equivalent to $C(X) \rtimes (\mathbb{R}^{n} \rtimes \Gamma)$ for some compact Hausdorff space $X$. This is the analog of the Cuntz-Li duality theorem for the algebra $\mathcal{U}_{\Gamma}$.

A matrix $A \in M_{d}(\mathbb{Z})$ is called an integer dilation matrix if all its eigenvalues are of absolute value greater than $1$. In \cite{Exel-Rae}, a purely infinite simple $C^{*}$-algebra associated to an integer dilation matrix was studied and its $K$-groups were computed.  Their computation depends on realising the $C^{*}$-algebra as a Cuntz-Pimsner algebra and by a careful examination of the six term sequence coming from its Toeplitz extension.
In \cite{Rae-Laca-Rammage}, a presentation of this algebra was obtained in terms of generators and relations.
For the group $\Gamma:=\{(A^{t})^{r}: r \in \mathbb{Z}\} \cong \mathbb{Z}$, denote the Cuntz-Li algebra $\mathcal{U}_{\Gamma}$ by $\mathcal{U}_{A^{t}}$. The presentation given in \cite{Rae-Laca-Rammage} tells us that the $C^{*}$-algebra studied in \cite{Exel-Rae} is the Cuntz-Li algebra $\mathcal{U}_{A^{t}}$.

The purpose of this paper is to understand the $K$-groups of $\mathcal{U}_{A^{t}}$ in view of the Cuntz-Li duality theorem. The Cuntz-Li duality theorem in this case says that $\mathcal{U}_{A^{t}}$ is Morita-equivalent to a crossed product algebra $(C(X) \rtimes \mathbb{R}^{d}) \rtimes \mathbb{Z}$ for some compact Hausdorff space. We compute the $K$-groups using the Pimsner-Voiculescu sequence. I believe that this computation will be of independent interest for the following two reasons:
\begin{enumerate}
 \item The $K$-groups of $\mathcal{U}_{A^{t}}$ depends on both $d$ and $\sgn(\det(A))$. (Cf. \cite{Exel-Rae}). The dependence on $d$ is due to the Connes-Thom isomorphism between $K_{*}(C(X))$ and $K_{*}(C(X) \rtimes \mathbb{R}^{d})$. Also the Connes-Thom map commutes with the action of $\mathbb{Z}$ if $\sgn(\det(A))>0$ and anticommutes if $\sgn(\det(A))$. This explains the dependence on $\sgn(\det(A))$.
\item It is mentioned in the introduction of \cite{Cuntz-Li-1} that the duality theorem enables one to use homotopy type arguments which makes it possible to compute the $K$-groups. We see the same kind of phenomenon here as well. (Cf. Lemma \ref{Trivial case}.)
 \end{enumerate}

Let $A \in GL_{n}(\mathbb{R})$. Consider the semidirect product $\mathbb{R}^{n} \rtimes \mathbb{Z}$ where $\mathbb{Z}$ acts on $\mathbb{R}^{n}$ by matrix multiplication by $A$. Let $B$ be a $C^{*}$-algebra on which $\mathbb{R}^{n} \rtimes \mathbb{Z}$. The crossed product $B \rtimes (\mathbb{R}^{n} \rtimes \mathbb{Z})$ is isomorphic to $(B \rtimes \mathbb{R}^{n}) \rtimes \mathbb{Z}$. In section 2 and 3, we write down the Pimsner-Voiculescu sequence for $(B \rtimes \mathbb{R}^{n}) \rtimes \mathbb{Z}$ after identifying the crossed product $B \rtimes \mathbb{R}^{n}$ with $B$ upto $KK$-equivalence. We show that the Connes-Thom isomorphism commutes with the action of $\mathbb{Z}$ if $\det(A)>0$ and anticommutes if $\det(A)<0$. In section 4 and 5, the $K$-groups of $\mathcal{U}_{A^{t}}$ were (re)computed.

\section{Preliminaries }
We use this section to fix notations and recall a few preliminaries.

Let $A \in M_{n}(\mathbb{R})$ be such that $\det(A) \neq 0$. We think of elements of $\mathbb{R}^{n}$ as column vectors. Thus the matrix $A$ induces an action of $\mathbb{Z}$ on $\mathbb{R}^{n}$ by left multiplication. The generator $1 \in \mathbb{Z}$ acts on $\mathbb{R}^{n}$ by $1.v=Av$ for $v \in \mathbb{R}^{n}$. Consider the semidirect product $\mathbb{R}^{n} \rtimes \mathbb{Z}$. 

Let $B$ be a $C^{*}$-algebra. A strongly continuous action of $\mathbb{R}^{n} \rtimes \mathbb{Z}$ on $B$ is equivalent to providing a pair  $(\alpha,\tau)$ where $\alpha$ is a strongly continuous action of $\mathbb{R}^{n}$ on $B$ and $\tau$ is an automorphism of $B$ such that 
\[\tau \alpha_{\xi}=\alpha_{A\xi}\tau\] for every $\xi \in \mathbb{R}^{n}$. If $(\alpha,\tau)$ is such a pair, we write $\alpha \rtimes \tau$ for the action of $\mathbb{R}^{n} \rtimes \mathbb{Z}$. Also the automorphism $\tau$ induces an action, denoted $\tilde{\tau}$, on the crossed product $B \rtimes_{\alpha} \mathbb{R}^{n}$ given by 
\begin{align*}
 \widetilde{\tau}(b):&=\tau(b) \textrm{~if~} b \in B \\
  \widetilde{\tau}(U_{\xi}):&=U_{A\xi} \textrm{~for~} \xi \in \mathbb{R}^{n}
\end{align*}
where $U_{\xi}$ denotes the canonical unitary in $\mathcal{M}(B \rtimes_{\alpha} \mathbb{R}^{n})$. Moreover the crossed product $B \rtimes_{\alpha \rtimes \tau} (\mathbb{R}^{n} \rtimes \mathbb{Z})$ is isomorphic to $(B \rtimes_{\alpha} \mathbb{R}^{n}) \rtimes_{\widetilde{\tau}} \mathbb{Z}$.

The Pimsner-Voiculescu sequence gives the following six-term exact sequence.

\begin{equation*}
\label{indexmap}
\def\labelstyle{\scriptstyle}
\xymatrix@C=25pt@R=20pt{
K_0(B \rtimes_{\alpha} \mathbb{R}^{n})\ar[r]^{1-\widetilde{\tau}_{*}}& K_0(B \rtimes_{\alpha} \mathbb{R}^{n})\ar[r]& K_0(B \rtimes_{\alpha \rtimes \tau}(\mathbb{R}^{n} \rtimes \mathbb{Z})) \ar[d] \\
K_1(B \rtimes_{\alpha \rtimes \tau}(\mathbb{R}^{n} \rtimes \mathbb{Z}))\ar[u]& K_1(B \rtimes_{\alpha} \mathbb{R}^{n} )\ar[l] & K_1(B \rtimes_{\alpha} \mathbb{R}^{n})\ar[l]^{1-\widetilde{\tau}_{*}} 
}.
\end{equation*}

But by the Connes-Thom isomorphism, we can replace $K_{i}(B \rtimes_{\alpha} \mathbb{R}^{n})$ by $K_{i+n}(B)$. Let $C_{n,i}:K_{i}(B) \to K_{i+n}(B\rtimes_{\alpha} \mathbb{R}^{n})$ be the Connes-Thom map. Now we can state our main theorem.

\begin{thm}
\label{main theorem}
 For $i=1,2$, $C_{n,i}\tau_{*}~=~\epsilon \widetilde{\tau}_{*}C_{n,i}$ where $\epsilon$ is given by 
\begin{displaymath}
\begin{array}{lll}
\epsilon&:=\left\{\begin{array}{lll}
                                
1 & \mbox{~if~} & \det(A) > 0 \\
                                 -1  & \mbox{~if~} & \det(A)<0 
                                 \end{array} \right. 
\end{array}
\end{displaymath}
\end{thm}

The following is an immediate corollary to Theorem \ref{main theorem}.

\begin{crlre}
\label{Main corollary}
 Let $(\alpha,\tau)$ be a strongly continuous action of $\mathbb{R}^{n} \rtimes \mathbb{Z}$ on a $C^{*}$-algebra $B$. Then there exists a six term exact sequence 
\begin{equation*}
\label{indexmap}
\def\labelstyle{\scriptstyle}
\xymatrix@C=25pt@R=20pt{
K_n(B)\ar[r]^{1-\epsilon \tau_{*}}& K_n(B)\ar[r]& K_0(B \rtimes_{\alpha \rtimes \tau}(\mathbb{R}^{n} \rtimes \mathbb{Z})) \ar[d] \\
K_1(B \rtimes_{\alpha \rtimes \tau}(\mathbb{R}^{n} \rtimes \mathbb{Z}))\ar[u]& K_{n+1}(B  )\ar[l] & K_{n+1}(B)\ar[l]^{1-\epsilon \tau_{*}} 
}.
\end{equation*}
where $\epsilon=\sgn(\det(A))$.

\end{crlre}

\section{Proof of Theorem \ref{main theorem}}
We use $KK$-theory to prove this. All our algebras are ungraded.
We denote the interior Kasparov product \[KK^{(i)}(A,B) \times KK^{(j)}(B,C) \to KK^{(i+j)}(A,C)\] by $\sharp$ and the external Kasparov product \[KK^{(i)}(A_{1},A_{2}) \times KK^{(j)}(B_{1},B_{2}) \to KK^{(i+j)}(A_{1}\otimes A_{2},B_{1} \otimes B_{2})\] by $\widehat{\otimes}$. We will also identify $K_{i}(B)$ with $KK^{(i)}(\mathbb{C},B)$.
Also if $\phi:B_{1} \to B_{2}$ is a $C^{*}$-algebra homomorphism then we denote the $KK$-element $(B_{2},\phi,0)$ in $KK^{(0)}(B_{1},B_{2})$ by $[\phi]$.

Under this identification, the Connes-Thom isomorphism is given by $C_{n}(x)=x ~\sharp~ t_{\alpha}$ where $t_{\alpha} \in KK^{n}(B,B \rtimes_{\alpha} \mathbb{R}^{n})$ is the Thom element.

Now it is immediate that Theorem \ref{main theorem} is equivalent to the following theorem.

\begin{thm}
 \label{Main theorem}
One has $[\tau]~\sharp~t_{\alpha}= \epsilon~ t_{\alpha} ~\sharp~ [\widetilde{\tau}]$ where $\epsilon=\sgn(\det(A))$.
\end{thm}

A bit of notation. If $X \in GL_{n}(\mathbb{R})$ then $X$ induces an automorphism  $\phi_{X}$ on $C_{0}\mathbb{R}^{n})$ given by $(\phi_{X}f)(v):=f(Xv)$. Let $b_{n} \in K_{n}(C_{0}(\mathbb{R}^{n}))$ be the Bott element. We denote the image $\phi_{X*}(b_{n}) \in K_{n}(C_{0}(\mathbb{R}^{n}))$ simply by $X_{*}(b_{n})$.

First let us dispose of the case when the action of $\mathbb{R}^{n}$ is trivial. For the trivial action the crossed product  $B \rtimes_{\alpha} \mathbb{R}^{n}$ is isomorphic to $B \otimes C_{0}(\mathbb{R}^{n})$ and $t_{\text{trivial}}=1_{B} \widehat{\otimes} b_{n}$.

\begin{lmma}
\label{Trivial case}
 If the action of $\mathbb{R}^{n}$ is trivial, then  $[\tau] ~\sharp~ t_{\text{trivial}}= \epsilon~ t_{\text{trivial}}~ \sharp~ [\widetilde{\tau}]$ where $\epsilon=\sgn(\det(A))$.
\end{lmma}
\textit{Proof.} Note that \begin{align*}
[\tau]~\sharp~t_{\text{trivial}}&=[\tau]~ \widehat{\otimes}~ b_{n} \\                           
t_{\text{trivial}}~ \sharp~ [\widetilde{\tau}]&= [\tau] ~\widehat{\otimes}~ A^{t}_{*}(b_{n})
                          \end{align*}
Thus we only need to prove that $A_{*}^{t}(b_{n})= \epsilon b_{n}$ where $\epsilon=\sgn(\det(A))$.

If $\det(A)>0$, then $A^{t}$ is homotopic to identity in $GL_{n}(\mathbb{R})$. Hence $A^{t}_{*}(b_{n})=b_{n}$. 

If $\det(A)<0$, then $A^{t}$ is homotopic to $\begin{pmatrix}
              -1 & 0 \\
                0 & Id_{n-1}                                                                                                                                                
                                 \end{pmatrix}$ in $GL_{n}(\mathbb{R})$. The Bott element $b_{n}=b_{1} \widehat{\otimes} b_{1} \widehat{\otimes} \cdots \widehat{\otimes} b_{1}$, it follows that the matrix $\begin{pmatrix}
              -1 & 0 \\
                0 & Id_{n-1}                                                                                                                                                
                                 \end{pmatrix}$ sends $b_{n}$ to $-b_{n}$. As a consequence, we have $A^{t}_{*}(b_{n})=-b_{n}$ if $\det(A)<0$.
This completes the proof. \hfill $\Box$

Now by a homotopy argument, the argument that is used in Theorem 2 of \cite{Fack-Skandalis}, we reduce Theorem \ref{Main theorem} to Lemma \ref{Trivial case}.

For $s \in [0,1]$, let $\alpha^{s}$ be the action of $\mathbb{R}^{n}$ on $B$ defined by $\alpha^{s}_{\xi}(b):=\alpha_{s\xi}(b)$. Note that $\alpha^{1}=\alpha$ and $\alpha^{0}$ gives the trivial action. Observe that $\tau \alpha^{s}_{\xi}=\alpha^{s}_{A\xi}\tau$. For $s \in [0,1]$, denote the automorphism $\tau$ by $\tau^{s}$ and the automorphism induced by $\tau$ on $B \rtimes_{\alpha^{s}}\mathbb{R}^{n}$ by $\widetilde{\tau}^{s}$.

Let $IB:=C[0,1]\otimes B$. Consider the action $\underline{\alpha}$ of $\mathbb{R}^{n}$ and the automorphism $\underline{\tau}$ on $IB$ defined by 
\begin{align*}
 \underline{\alpha}_{\xi}(f)(s):&=\alpha^{s}_{\xi}(f(s))\\
\underline{\tau}(f)(s)&=\tau(f(s))
\end{align*}

Observe that for $\xi \in \mathbb{R}^{n}$, $\underline{\tau}~\underline{\alpha}_{\xi}=\underline{\alpha}_{A\xi}~\underline{\tau}$. The automorphism $\underline{\tau}$ induces an automorphism on $IB \rtimes_{\underline{\alpha}} \mathbb{R}^{n}$ and we denote it by $\widetilde{\underline{\tau}}$.

For $s \in [0,1]$, let $\epsilon_{s}:IB \to B$ be the evaluation map. Then $\epsilon_{s}:(IB,\underline{\alpha}) \to (B,\alpha^{s})$ is equivariant. We denote the induced map from $IB \rtimes_{\underline{\alpha}} \mathbb{R}^{n}$ to $B \rtimes_{\alpha^{s}} \mathbb{R}^{n}$ by $\widehat{\epsilon_{s}}$.

Also for $s \in [0,1]$, $\widehat{\epsilon_{s}} \circ \widetilde{\underline{\tau}}=\widetilde{\tau}^{s}\circ \widehat{\epsilon_{s}}$.

\begin{lmma}
 For $s \in [0,1]$, the element $[\widehat{\epsilon}_{s}] \in KK^{(0)}(IB \rtimes_{\underline{\alpha}}\mathbb{R}^{n},B \rtimes_{\alpha^{s}} \mathbb{R}^{n})$ is a $KK$-equivalence.
\end{lmma}
\textit{Proof.} Observe that \[
                               t_{\underline{\alpha}} ~\sharp~ [\widehat{\epsilon_{s}}]= [\epsilon_{s}] ~\sharp~[t_{\alpha^{s}}].
                             \]
Since $[\epsilon_{s}] \in KK^{(0)}(IB,B)$ and the Thom elements   are $KK$-equivalences, it follows that $[\widehat{\epsilon_{s}}]$ is a $KK$-equivalence. This completes the proof. \hfill $\Box$

\begin{ppsn}
 \label{homotopy argument}
The following are equivalent. Recall that $\epsilon=\sgn(\det(A))$.
\begin{enumerate}
 \item For every $s \in [0,1]$, $[\tau^{s}]~\sharp~t_{\alpha^{s}}=\epsilon~t_{\alpha^{s}}~\sharp~[\widetilde{\tau}^{s}]$.
 \item There exists $s \in [0,1]$ such that $[\tau^{s}]~\sharp~t_{\alpha^{s}}=\epsilon~t_{\alpha^{s}}~\sharp~[\widetilde{\tau}^{s}]$.
\item The Kasparov product $[\underline{\tau}]~\sharp~[t_{\underline{\alpha}}]=\epsilon~t_{\underline{\alpha}}~\sharp~[\widetilde{\underline{\tau}}]$.
\end{enumerate}
\end{ppsn}
\textit{Proof.}  Let $s \in [0,1]$ be given. Observe the following.
\begin{align*}
[\underline{\tau}]~\sharp~t_{\underline{\alpha}}=&\epsilon~t_{\underline{\alpha}}~\sharp~[\widetilde{\underline{\tau}}] \\
\Leftrightarrow & [\underline{\tau}]~\sharp~ t_{\underline{\alpha}}~\sharp~[\widehat{\epsilon_{s}}]=\epsilon~t_{\underline{\alpha}}~\sharp~[\widetilde{\underline{\tau}}]~\sharp~[\widehat{\epsilon_{s}}] \text{~(Since $[\widehat{\epsilon_{s}}]$ is a $KK$-equivalence.)} \\
\Leftrightarrow & [\underline{\tau}]~\sharp~[\epsilon_{s}]~\sharp~t_{\alpha^{s}}=\epsilon~t_{\underline{\alpha}}~\sharp~[\widehat{\epsilon_{s}} \circ \widetilde{\underline{\tau}}] \\
\Leftrightarrow &
[\epsilon_{s} \circ \underline{\tau}]~\sharp~t_{\alpha^{s}}=\epsilon~t_{\underline{\alpha}}~\sharp~[\widetilde{\tau^{s}} \circ \widehat{\epsilon_{s}}] \\
\Leftrightarrow &
[\tau^{s} \circ \epsilon_{s}]~\sharp~t_{\alpha^{s}}= \epsilon~t_{\underline{\alpha}}~\sharp~[\widehat{\epsilon_{s}}]~\sharp~[\widetilde{\tau^{s}}] \\
\Leftrightarrow &
[\epsilon_{s}]~\sharp~[\tau^{s}]~\sharp t_{\alpha^{s}}=\epsilon~[\epsilon_{s}]~\sharp~t_{\alpha^{s}}~\sharp~[\widetilde{\tau^{s}}] \\
\Leftrightarrow &
[\tau^{s}]~\sharp~[t_{\alpha^{s}}]=\epsilon~t_{\alpha^{s}}~\sharp~[\widetilde{\tau^{s}}] \text{~(Since $[\epsilon_{s}]$ is a $KK$-equivalence.)}
\end{align*}
The proof is now complete. \hfill $\Box$.

Now Theorem \ref{Main theorem} follows from Proposition \ref{homotopy argument} and Lemma \ref{Trivial case}.

\section{The Cuntz-Li algebra associated to an integer dilation matrix}

As an application of Corollary \ref{Main corollary}, we recompute the $K$-theory of the $C^{*}$-algebra associated to an integer dilation matrix, studied in \cite{Exel-Rae}. Let us recall the $C^{*}$-algebra considered in \cite{Exel-Rae}.

 Let $A \in M_{d}(\mathbb{Z})$ be an integer dilation matrix i.e. all the eigen values of $A$ are of absolute value greater than $1$. The matrix $A$ acts on $\mathbb{R}^{d}$ by matrix multiplication and leaves $\mathbb{Z}^{d}$ invariant. Denote the resulting endomorphism on $\mathbb{T}^{d}:=\mathbb{R}^{d}/\mathbb{Z}^{d}$ by $\sigma_{A}$. The map $\sigma_{A}$ is a surjective and has finite fibres. Denote the map $C(\mathbb{T}^{d}) \ni f \to \in f\circ \sigma_{A} C(\mathbb{T}^{d})$ by $\alpha_{A}$. Consider the transfer operator $L:C(\mathbb{T}^{d}) \to C(\mathbb{T}^{d})$ defined by 
\[
 L(f)(x):=\frac{1}{|\sigma_{A}^{-1}(x)|}\sum_{\sigma_{A}(y)=x}f(y)
\]
Then $L$ satisfies the condition $L(\alpha_{A}(f)g)=fL(g)$ for $f,g \in C(\mathbb{T}^{d})$. In \cite{Exel-Rae}, the Exel Crossed product $C(\mathbb{T}^{d}) \rtimes_{\alpha_{A},L} \mathbb{N}$ was viewed as a Cuntz-Pimsner algebra $\mathcal{O}(M_{L})$ of a suitable Hilbert $C(\mathbb{T}^{d})$ bimodule $M_{L}$.

By a careful examination of the six term sequence (and the maps involved) associated to the  exact sequence $0 \to Ker(Q) \to \mathcal{T}(M_{L}) \to \mathcal{O}(M_{L})$, the $K$-groups of $C(\mathbb{T}^{d}) \rtimes_{\alpha_{A}} \mathbb{N}$ were computed in \cite{Exel-Rae}.

For our purposes, the following description of $\mathcal{O}(M_{L})$ in terms of generators and relations is more relevant. Let us recall the following proposition from \cite{Rae-Laca-Rammage} (Proposition 3.3, Page 6 ).

\begin{ppsn}[\cite{Rae-Laca-Rammage}]
\label{Laca-Raeburn-Rammage}
The Exel's crossed product $C(\mathbb{T}^{d}) \rtimes_{\alpha_{A},L} \mathbb{N}$ is the universal $C^{*}$-algebra generated by an isometry $v$ and unitaries $\{u_{m}:m\in \mathbb{Z}^{d} \}$ satisfying  the following relations. 
\begin{align*}
 u_{m}u_{n}&=u_{m+n} \\
 vu_{m}&=u_{A^{t}m}v \\
 \sum_{ m \in \Sigma} u_{m}vv^{*}u_{m}^{-1}&=1
\end{align*}
Here $\Sigma$ denotes a set of distinct coset representatives of the group $\mathbb{Z}^{d}/A^{t}\mathbb{Z}^{d}$.
\end{ppsn}

\begin{rmrk}
 The above relations are called condtions $(E1)$ and $(E3)$ in \cite{Rae-Laca-Rammage}. Condition $(E2)$ in \cite{Rae-Laca-Rammage} is implied by $(E1)$ and $(E3)$. 

For if $ m \notin A^{t}\mathbb{Z}^{d}$, the projections $vv^{*}$ and $u_{m}vv^{*}u_{m}^{-1}$ are orthogonal by $(E3)$. Hence $v^{*}u_{m}v=0$ if $m \notin A^{t}\mathbb{Z}^{d}$. If $ m \in A^{t}\mathbb{Z}^{d}$, then using $vu_{m}=u_{A^{t}m}v$, one obtains $v^{*}u_{m}v=u_{(A^{t})^{-1}m}$. Thus $(E1)$ and $(E3)$ implies $(E2)$.
\end{rmrk}

The following setup was initially considered in \cite{Quigg-Landstad}. Consider a semi-direct product $N \rtimes H$ and let $M$ be a normal subgroup. Let $P:=\{a \in H:aMa^{-1} \subset M\}$. Then $P$ is a semigroup containing the identity $e$. For $a \in P$, let $M_{a}=aMa^{-1}$.  Assume that the following holds.
\begin{enumerate}
 \item[(C1)] The group $H=PP^{-1}=P^{-1}P$.
  \item[(C2)] For every $a \in P$, the subgroup $aMa^{-1}$ is of finite index in $M$.
  \item[(C3)] The intersection $\displaystyle \bigcap_{a \in P}aMa^{-1} = \{e\}$ where $e$ denotes the identity element of $G$.
\end{enumerate}

\begin{dfn}
The Cuntz-Li algebra associated to the pair $(N \rtimes H,M)$  is the  the universal $C^{*}$-algebra generated by a set of isometries $\{s_{a}:a \in P\}$ and a set of unitaries $\{u(m): m \in M\}$ satisfying the following relations.
\begin{equation*}
 \begin{split}
  s_{a}s_{b}&=s_{ab}\\
  u(m)u(n)&=u(mn)\\
  s_{a}u(m)&=u(ama^{-1})s_{a}\\
  \displaystyle \sum_{k \in M/M_{a}}u(k)e_{a}u(k)^{-1}&=1
 \end{split}
\end{equation*}
where $e_{a}$ denotes the final projection of $s_{a}$. We denote the Cuntz-Li algebra associated to the pair $(N \rtimes H, M)$ by $U[N \rtimes H,M]$.
\end{dfn}

Let $A \in M_{d}(\mathbb{Z})$ be a dilation matrix. Then $A$ acts on $\mathbb{Q}^{d}$ by left multiplication. Consider the semidirect product $\mathbb{Q}^{d} \rtimes \mathbb{Z}$ and the normal subgroup $\mathbb{Z}^{d}$ of $\mathbb{Q}^{d}$. For this pair $(\mathbb{Q}^{d} \rtimes \mathbb{Z},~\mathbb{Z}^{d})$, $P=\{A^{r}:r \geq 0\} \cong \mathbb{N} $. Moreover conditions $(C1)-(C3)$ are satisfied. (See Example 2.6, Page 3 in \cite{Sundar-Cuntz-Li}.)
Let us denote the Cuntz-Li algebra $U[\mathbb{Q}^{d} \rtimes \mathbb{Z}, \mathbb{Z}^{d}]$ simply by $\mathcal{U}_{A}$.

By using the presentation (Cf. Prop. \ref{Laca-Raeburn-Rammage}) of the Exel's Crossed product $C(\mathbb{T}^{d}) \rtimes_{\alpha_{A},L} \mathbb{N}$ given in terms of isometries and unitaries, it is easy to verify that $C(\mathbb{T}^{d}) \rtimes_{\alpha_{A},L} \mathbb{N}$ is isomorphic to $\mathcal{U}_{A^{t}}$. 

Let us recall the Cuntz-Li duality result proved in \cite{Sundar-Cuntz-Li}. The proof is really a step by step adaptation of the arguments used in \cite{Cuntz-Li-1}.

Let $\displaystyle N_{A}:= \bigcup_{r=0}^{\infty}A^{-r}\mathbb{Z}^{d}$. Then $N_{A}$ is a subgroup of $\mathbb{R}^{d}$. Let $\mathbb{Z}$ act on $\mathbb{R}^{d}$ by left multiplication by $A^{t}$ and consider the semidirect product $\mathbb{R}^{d} \rtimes \mathbb{Z}$. The semidirect product $\mathbb{R}^{d} \rtimes \mathbb{Z}$ acts on the $C^{*}$-algebra $C^{*}(N_{A})$. The action $\alpha$ of $\mathbb{R}^{d}$ and the automorphism $\tau$, corresponding to the action of $\mathbb{Z}$, are given by
\begin{align*}
 \alpha_{\xi}(\delta_{v}):&=e^{-2\pi i <\xi,v>}\delta_{v} \mbox{~for~} \xi \in \mathbb{R}^{d}\\
  \tau(\delta_{v}):&=\delta_{A^{-1}v}
\end{align*}
where $\{\delta_{v}:v \in N_{A}\}$ denotes the canonical unitaries in $C^{*}(N_{A})$ and $<,>$ denotes the usual inner product on $\mathbb{R}^{d}$. Note that $\tau \alpha_{\xi}=\alpha_{A^{t}\xi}\tau$ for $\xi \in \mathbb{R}^{d}$.

The following proposition was proved in \cite{Sundar-Cuntz-Li}. (Cf. Theorem 8.2 and Proposition 8.6 in \cite{Sundar-Cuntz-Li})

\begin{ppsn}
\label{duality}
The $C^{*}$-algebra $\mathcal{U}_{A^{t}}$ is Morita-equivalent to $C^{*}(N_{A}) \rtimes_{\alpha \rtimes \tau} (\mathbb{R}^{d} \rtimes_{A^{t}} \mathbb{Z})$. 
\end{ppsn}

Now using the Morita-equivalence in Proposition \ref{duality} and using the version of Pimsner-Voiculescu exact sequence established in Corollary \ref{Main corollary}, the $K$-groups of $\mathcal{U}_{A^{t}}$ can be computed.

\section{K-groups of the Cuntz-Li algebra $\mathcal{U}_{A^{t}}$}

Recall that $N_{A}:=\bigcup_{r=0}^{\infty}A^{-r}\mathbb{Z}^{d}$. Set $N_{A}^{(r)}:=A^{-r}\mathbb{Z}^{d}$. Then $\{N_{A}^{(r)}\}_{r=0}^{\infty}$ forms an increasing sequence of subgroups, each isomorphic to $\mathbb{Z}^{d}$ and $N_{A}=\bigcup_{r=0}^{\infty}N_{A}^{(r)}$. Thus $C^{*}(N_{A})$ is the inductive limit of $C^{*}(N_{A}^{(r)}) \cong C(\mathbb{T}^{d})$. Thus $K_{*}(C^{*}(N_{A}))$ can be computed as the inductive limit of the $K$-groups of $C^{*}(N_{A}^{(r)}) \cong C(\mathbb{T}^{d})$.

Let us first recall the $K$-theory of $C(\mathbb{T}^{d}) \cong C^{*}(\mathbb{Z}^{d})$. It is well known and can be proved by the Kunneth formula that as a $\mathbb{Z}_{2}$ graded ring, $K_{*}(C^{*}(\mathbb{Z}^{d}))$ is isomorphic to the exterior algebra $ \Lambda^{*}(\mathbb{Z}^{d})$. The map $ \mathbb{Z}^{d} \ni e_{i} \to \delta_{e_{i}} \in K_{1}(C^{*}(\mathbb{Z}^{d}))$ extends to a graded ring isomorphism from $\Lambda^{*}(\mathbb{Z}^{d})$ to $K_{*}(C^{*}(\mathbb{Z}^{d}))$. 

\begin{rmrk}
The isomorphism $K_{*}(C(\mathbb{T}^{d})) \cong \Lambda^{*}\mathbb{Z}^{d}$  was also used in \cite{Exel-Rae}.
\end{rmrk}

Let us now fix some notations. For $0 \leq n \leq d$, let $A_{n}$ be the map on $\Lambda^{n}(\mathbb{Z}^{d})$ induced by $A$.
Thus $A_{0}=1$, $A_{1}=A$ and $A_{d}=\det(A)$. For a subset  $I$ of $\{1,2,\cdots,d\}$,  of cardinality $n$, (arranged in increasing order ), $I =\{i_{1}<i_{2}< \cdots, <i_{n} \}$, let $e_{I}:=e_{i_{1}}\wedge e_{i_{2}} \wedge \cdots \wedge e_{i_{n}}$. Then  $\{e_{I}:|I|=n\}$ is a basis for $\Lambda^{n}(\mathbb{Z}^{d})$. For subsets $J,K$ of $\{1,2,\cdots,d\}$ of size $n$, let $A_{J,K}$ be the submatrix of $A$ obtained by considering the rows coming from $J$ and columns coming from $K$. With respect to the basis $\{e_{I}:|I|=n\}$, the $(J,K)^{th}$ entry of the matrix corresponding to $A_{n}$ is $\det(A_{J,K})$. 

Note that for $n \geq 1$, $A_{n}$ is again an integer dilation matrix. For if we upper triangualise $A$, then w.r.t the basis $\{e_{I}:|I|=n\}$, arranged in lexicographic order, $A_{n}$ is upper triangular and the eigen values of $A_{n}$ are product of eigen values of $A$. Thus the eigen values of $A_{n}$ are of absolute value greater than $1$.
This fact was  used in \cite{Exel-Rae}. (Cf. Proposition 4.6. in \cite{Exel-Rae}.)

Let $n \in \{0,1,2,\cdots,n\}$. Consider $\Lambda^{n}(\mathbb{Z}^{d})$ as a subgroup of $\Lambda^{n}(\mathbb{Q}^{d})$. Then $A_{n}$ is invertible on $\Lambda^{n}(\mathbb{Q}^{d})$. Set $ \displaystyle \Gamma_{n}:=\bigcup_{r=0}^{\infty} A_{n}^{-r} (\Lambda^{n}(\mathbb{Z}^{d}))$.

\begin{ppsn}
\label{group algebra}
 The $K$-groups of the  $C^{*}$-algebra $C^{*}(N_{A})$ are given by 
\begin{align*}
  K_{0}(C^{*}(N_{A})) &\cong   \mathop{\bigoplus_{n \text{~even~}}}_{0 \leq n \leq d}\Gamma_{n}\\
  K_{1}(C^{*}(N_{A})) &\cong \bigoplus_{\substack {n \text{~odd~} \\ 0 \leq n \leq d}}\Gamma_{n} 
\end{align*}
\end{ppsn}
\textit{Proof.} Since $C^{*}(N_{A})$ is the inductive limit of $C^{*}(N_{A}^{(r)}) \cong C(\mathbb{T}^{d})$, it follows that $K_{*}(C^{*}(N_{A}))$ is the inductive limit of $K_{*}(C^{*}(N_{A}^{(r)})) \cong K_{*}(C(\mathbb{T}^{d}))$. Identify $K_{*}(C^{*}(N_{A}^{(r)}))$ with $\Lambda^{*}(\mathbb{Z}^{d})$ via the map $\delta_{A^{-r}(e_{i})} \to e_{i}$. With this identification the inclusion map $C^{*}(N_{A}^{(r)}) \to C^{*}(N_{A}^{-(r+1)})$ induces the map $\displaystyle \oplus_{ 1 \leq n \leq d}A_{n}$ at the $K$-theory level. (Reason: If we write $e_{j}$ as a linear combination of $\{A^{-1}e_{i}\}$ the matrix involved is just $A$. )

Thus we are left to show that the inductive limit of $(\bigoplus_{n} \Lambda^{n}\mathbb{Z}^{d},\bigoplus_{n}A_{n})_{r=0}^{\infty}$ is $\bigoplus_{n} \Gamma_{n}$. Again it is enough to show that the inductive limit of $(\Lambda^{n}\mathbb{Z}^{d},A_{n})_{r=0}^{\infty}$ is isomorphic to $\Gamma_{n}$. Let $H_{r}=\Lambda^{n}\mathbb{Z}^{d}$. If $v \in \Gamma_{n}$, write $v$ as $v=A_{n}^{-r}w$ with $w \in \Lambda^{n}(\mathbb{Z}^{d})$. The map $\Gamma_{n} \ni v \to w \in H_{r}$ is an isomorphism between $\Gamma_{n}$ and $ \displaystyle \lim_{r\to \infty}(\Lambda^{n}(\mathbb{Z}^{d}),A_{n})$. This completes the proof. \hfill $\Box$

Now let us calculate the automorphism $\tau$ on $C^{*}(N_{A})$. Recall that $\tau$ on the generating unitaries is given by $\tau(\delta_{v})=\delta_{A^{-1}v}$. Thus, it is immediate and not difficult to see that $\tau$ induces the map $\bigoplus_{n}A_{n}^{-1}$ on $\bigoplus_{n}\Gamma_{n}$ when one identifies $K_{*}(C^{*}(N_{A}))$ with $\bigoplus_{n}\Gamma_{n}$ ( together with their $\mathbb{Z}_{2}$ grading).

We need one more lemma.

\begin{lmma}
 Let $1 \leq n \leq d$. The natural map \[\frac{\Lambda^{n}(\mathbb{Z}^{d})}{(1-\epsilon A_{n})(\Lambda^{n}(\mathbb{Z}^{d}))} \to \frac{\Gamma_{n}}{(1-\epsilon A_{n}) \Gamma_{n}}\] is an isomorphism for $\epsilon \in \{1,-1\}$.
\end{lmma}
\textit{Proof.} Let us denote $\Lambda^{n}(\mathbb{Z}^{d})$ by $V_{n}$. We will give a proof only for $\epsilon =1$. The case $\epsilon=-1$ is similar and we leave its proof to the reader.

Observe that for $1-A_{n}^{r}=(1-A_{n})(1+A_{n}+A_{n}^{2}+\cdots +A_{n}^{r-1})$ for $r \geq 0$. Thus for every $ r \geq 0$, there exists a polynomial $p_{r}(x)$ with integer co-efficients such that $p_{r}(A_{n})(1-A_{n})+A_{n}^{r}=1$.

\textit{Surjectivity:} Let $v \in \Gamma_{n}$ be given. By definition, there exists $r \geq 0$ and $w \in V_{n}$ such that $v=A_{n}^{-r}w$. Hence \begin{align*}
v&=A_{n}^{-r}(A_{n}^{r}+p_{r}(A_{n})(1-A_{n}))w \\                                                                                                                                                   
 &=w+ (1-A_{n})p_{r}(A_{n})A_{n}^{-r}w                                                                                                                                                \end{align*}
Hence $v \equiv w \mod (1-A_{n})\Gamma_{n}$. This proves the surjectivity of the given map.

\textit{Injectivity:} Let $v \in V_{n}$ be such that $v \equiv 0 \mod (1-A_{n})\Gamma_{n}$. This implies that there exists $r \geq 0$ and $w \in \Lambda^{n}(\mathbb{Z}^{d})$ such that $v=(1-A_{n})A_{n}^{-r}w$. Hence $(1-A_{n})w=A_{n}^{r}v$. Now observe that 
\begin{align*}
 w &= (p_{r}(A_{n}(1-A_{n})+A_{n}^{r})w \\
   &=p_{r}(A_{n})A_{n}^{r}v+A_{n}^{r}w \\
   &=A_{n}^{r}(p_{r}(A_{n})v+w)
\end{align*}
Hence $A_{n}^{-r}w=p_{r}(A_{n})v+w \in V_{n}$. Hence $v \equiv 0 \mod (1-A_{n})V_{n}$. This proves the injectivity part. This completes the proof. \hfill $\Box$

We denote both the abelian groups  $\frac{\Lambda^{n}(\mathbb{Z}^{d})}{(1-\epsilon A_{n})(\Lambda^{n}(\mathbb{Z}^{d}))}$ and $\frac{\Gamma_{n}}{(1-\epsilon A_{n}) \Gamma_{n}}$ by $coker(1-\epsilon A_{n})$.

\begin{thm}
\label{Cuntz-Li-Ktheory}
 Let $A \in M_{d}(\mathbb{Z})$ be an integer dilation matrix. The $K$-groups of the Cuntz-Li algebra~ $\mathcal{U}_{A^{t}}$ are as follows. 

\begin{enumerate}
\item If $d$ is even and $\det(A)>0$, then 
\begin{align*}
 K_{0}(\mathcal{U}_{A^{t}}) & \displaystyle ~\cong ~\bigoplus_{\substack{ n \text{~ even~} \\ 0 \leq n \leq d}}coker(1-A_{n}), \textrm{~and~} \\
K_{1}(\mathcal{U}_{A^{t}}) & ~\cong~ \bigoplus_{\substack{n \text{~odd~} \\ 0 \leq n \leq d}}coker(1-A_{n}) \oplus \mathbb{Z}.
\end{align*}
\item If $d$ is even and $\det(A)<0$ then
\begin{align*}
 K_{0}(\mathcal{U}_{A^{t}}) & \displaystyle ~\cong ~\bigoplus_{\substack{ n \text{~ even~} \\ 0 \leq n \leq d}}coker(1+A_{n}), \textrm{~and~} \\
K_{1}(\mathcal{U}_{A^{t}}) & ~\cong~ \bigoplus_{\substack{n \text{~odd~} \\ 0 \leq n \leq d}}coker(1+A_{n}). 
\end{align*}
\item If $d$ is odd and $\det(A)>0$ then 
\begin{align*}
 K_{0}(\mathcal{U}_{A^{t}}) & \displaystyle ~\cong ~\bigoplus_{\substack{ n \text{~ odd~} \\ 0 \leq n \leq d}}coker(1-A_{n}) \oplus \mathbb{Z}, \textrm{~and~}\\
K_{1}(\mathcal{U}_{A^{t}}) & ~\cong~ \bigoplus_{\substack{n \text{~even~} \\ 0 \leq n \leq d}}coker(1-A_{n}). 
\end{align*}
\item If $d$ is odd and $\det(A)<0$ then
\begin{align*}
 K_{0}(\mathcal{U}_{A^{t}}) & \displaystyle ~\cong ~\bigoplus_{\substack{ n \text{~ odd~} \\ 0 \leq n \leq d}}coker(1+A_{n}), \textrm{~and~} \\
K_{1}(\mathcal{U}_{A^{t}}) & ~\cong~ \bigoplus_{\substack{n \text{~even~} \\ 0 \leq n \leq d}}coker(1+A_{n}). 
\end{align*}
\end{enumerate}
\end{thm}
\textit{Proof.} Our main tool is Corollary \ref{Main corollary} and the Morita equivalence between $\mathcal{U}_{A^{t}}$ and $C^{*}(N_{A}) \rtimes (\mathbb{R}^{n} \rtimes \mathbb{Z})$.

\textbf{If $d$ is even and $\det(A)>0$} then by Corollary \ref{Main corollary} , one has the following six term exact sequence.

\begin{equation*}
\label{indexmap}
\def\labelstyle{\scriptstyle}
\xymatrix@C=25pt@R=20pt{
K_0(C^{*}(N_{A}))\ar[r]^{1- \tau_{*}}& K_0(C^{*}(N_{A}))\ar[r]& K_0(\mathcal{U}_{A^{t}}) \ar[d] \\
K_1(\mathcal{U}_{A^{t}})\ar[u]& K_{1}(C^{*}(N_{A})) \ar[l] & K_{1}(C^{*}(N_{A}))\ar[l]^{1- \tau_{*}} 
}.
\end{equation*}
Now by Prop. \ref{group algebra}, the above six term sequence becomes

\begin{equation*}
\label{indexmap}
\def\labelstyle{\scriptstyle}
\xymatrix@C=45pt@R=35pt{
\displaystyle \bigoplus_{\substack{n \text{~even~} \\ 0 \leq n \leq d}}\Gamma_{n}\ar[r]^{ 1- \oplus_{n}A_{n}^{-1}}&  \displaystyle \bigoplus_{\substack{n \text{~even~} \\ 0 \leq n \leq d}}\Gamma_{n} \ar[r] & K_0(\mathcal{U}_{A^{t}}) \ar[d] \\
K_1(\mathcal{U}_{A^{t}})\ar[u]& \displaystyle \bigoplus_{\substack{n \text{~odd~} \\ 0 \leq n \leq d}}\Gamma_{n} \ar[l] & \displaystyle \bigoplus_{\substack{n \text{~odd~} \\ 0 \leq n \leq d}}\Gamma_{n}\ar[l]^{1- \oplus_{n}A_{n}^{-1}} 
}.
\end{equation*}
Now for $n \geq 1$, $A_{n}$ is a dilation matrix and thus $ker(1-A_{n}^{-1})=0$ if $n \geq 1$. Hence we conclude from the above six term sequence that \[K_{0}(\mathcal{U}_{A^{t}}) \equiv \bigoplus_{\substack{n \text{~even~} \\ 0 \leq n \leq d}}coker(1-A_{n}^{-1}).\]
Since $A_{n}$ is invertible and $(1-A_{n})=-A_{n}(1-A_{n}^{-1})$, it follows that $coker(1-A_{n}) \equiv coker(1-A_{n}^{-1})$. Thus \[K_{0}(\mathcal{U}_{A^{t}}) \equiv \bigoplus_{\substack{n \text{~even~} \\ 0 \leq n \leq d}}coker(1-A_{n}^{-1})  \equiv \bigoplus_{\substack{n \text{~even~} \\ 0 \leq n \leq d}}coker(1-A_{n}) .\]

Now $A_{0}=1$ and hence $\displaystyle \bigoplus_{n \text{~even}}ker(1-A_{n}^{-1})=\mathbb{Z}$. Again the six term sequence gives the following short exact sequence.
\[
 0 \longrightarrow \displaystyle \bigoplus_{\substack{n \text{~odd~} \\ 0 \leq n \leq d}}coker(1-A_{n}^{-1}) \longrightarrow K_{1}(\mathcal{U}_{A^{t}}) \longrightarrow \mathbb{Z}
\]
Since $\mathbb{Z}$ is free, it follows that \[
                                             K_{1}(\mathcal{U}_{A^{t}}) \equiv \displaystyle \bigoplus_{\substack{n \text{~odd~} \\ 0 \leq n \leq d}}coker(1-A_{n}^{-1}) \oplus \mathbb{Z} 
                                                 \equiv \displaystyle \bigoplus_{\substack{n \text{~odd~} \\ 0 \leq n \leq d}}coker(1-A_{n}) \oplus \mathbb{Z}.
                                                     \]

\textbf{If $d$ is even and $\det(A)<0$} then by Corollary \ref{Main corollary} and \ref{group algebra}, we get  the following six term sequence.

\begin{equation*}
\label{indexmap}
\def\labelstyle{\scriptstyle}
\xymatrix@C=45pt@R=35pt{
\displaystyle \bigoplus_{\substack{n \text{~even~} \\ 0 \leq n \leq d}}\Gamma_{n}\ar[r]^{1 + \oplus_{n}A_{n}^{-1}}&  \displaystyle \bigoplus_{\substack{n \text{~even~} \\ 0 \leq n \leq d}}\Gamma_{n} \ar[r] & K_0(\mathcal{U}_{A^{t}}) \ar[d] \\
K_1(\mathcal{U}_{A^{t}})\ar[u]& \displaystyle \bigoplus_{\substack{n \text{~odd~} \\ 0 \leq n \leq d}}\Gamma_{n} \ar[l] & \displaystyle \bigoplus_{\substack{n \text{~odd~} \\ 0 \leq n \leq d}}\Gamma_{n}\ar[l]^{1+ \oplus_{n}A_{n}^{-1}} 
}.
\end{equation*}

Again for $n \geq 1$, $A_{n}$ is a dilation matrix and for $n=0$, $A_{0}=1$. Hence for $0 \leq n \leq d$, $ker(1+A_{n}^{-1})=0$.  Thus the above six term sequence implies that 
\begin{align*}
 K_{0}(\mathcal{U}_{A^{t}}) &\equiv \bigoplus_{\substack{n \text{~even~} \\ 0 \leq n \leq d}}coker(1+A_{n}^{-1}) \equiv \bigoplus_{\substack{n \text{~even~} \\ 0 \leq n \leq d}}coker(1+A_{n}) \\
 K_{1}(\mathcal{U}_{A^{t}}) &\equiv \bigoplus_{\substack{n \text{~odd~} \\ 0 \leq n \leq d}}coker(1+A_{n}^{-1}) \equiv \bigoplus_{\substack{n \text{~odd~} \\ 0 \leq n \leq d}}coker(1+A_{n})
\end{align*}

The case when $d$ is odd is similar (again an application of Corollary \ref{Main corollary}) and we leave the details to the reader. This completes the proof. \hfill $\Box$

The rest of this section is devoted to reconciling Theorem \ref{Cuntz-Li-Ktheory} with the result obtained in \cite{Exel-Rae}. More precisely with Theorem 4.9 of \cite{Exel-Rae}. Let us recall the notations as in \cite{Exel-Rae}.

For a subset $K=\{k_{1}<k_{2}<\cdots <k_{n}\}$ of $\{1,2,\cdots,d\}$, denote the complement arranged in increasing order by  $K^{'}$ and let $K^{'}=\{k_{n+1}<k_{n+2}<\cdots<k_{d}\}$. Denote the permutation $i \to k_{i}$ by $\tau_{K}$. For a permutation $\sigma$, $\sgn(\sigma)$ is $1$ if $\sigma$ is even and $-1$ if $\sgn(\sigma)$ is odd.
Also recall that if $K$ and $J$ are subsets of size $n$, then $A_{K,J}$ is the matrix obtained from $A$ by considering the rows from $K$ and columns from $J$.

 For $0 \leq n \leq d$, let $\widetilde{B}_{n}$ be the $\binom{d}{n} \times \binom{d}{n}$ matrix defined as follows. ( We index the columns and rows by subsets of $\{1,2,\cdots,d\}$ of size $n$.) The $(K,L)^{th}$ entry of $\widetilde{B}_{n}$ is $\sgn(\tau_{K}\tau_{L})\det(A_{K^{'},L^{'}})$.

The matrices $B_{n}$ as defined in \cite{Exel-Rae} (Prop 4.6.) are then given by $B_{n}=\sgn(\det(A))\widetilde{B}_{n}$.  Denote the matrix whose $(K,L)^{th}$ entry is $det(A_{K^{'},L^{'}})$ by $C_{n}$. By convention, $\widetilde{B}_{d}=1=C_{d}$. Note that $\widetilde{B}_{n}$ and $C_{n}$ are conjugate over $\mathbb{Z}$. For the matrix $diag(\sgn(\tau_{K}))$ conjugates $\widetilde{B}_{n}$ to $C_{n}$.

Let $U_{n}:\Lambda^{n}(\mathbb{Z}^{d}) \to \Lambda^{d-n}(\mathbb{Z}^{d})$ be defined by $U_{n}e_{I}:=e_{I^{'}}$. Then $U_{n}$ is invertible and $U_{n}C_{n}U_{n}^{-1}=A_{d-n}$. Since $\sgn(\det(A))C_{n}$ is conjugate (over $\mathbb{Z}$) to ${B}_{n}$, it follows that $\sgn(\det(A))A_{d-n}$ is conjugate (over $\mathbb{Z}$) to $B_{n}$.

Now Theorem \ref{Cuntz-Li-Ktheory} can be restated, in terms of the matrices $B_{n}$'s, as in the following proposition. This is exactly Theorem 4.9 of \cite{Exel-Rae}.

\begin{thm}
\label{Cuntz-Li-Ktheory}
 Let $A \in M_{d}(\mathbb{Z})$ be an integer dilation matrix. The $K$-groups of the Cuntz-Li algebra~ $\mathcal{U}_{A^{t}}$ are as follows. 

\begin{enumerate}
\item If $d$ is even and $\det(A)>0$, then 
\begin{align*}
 K_{0}(\mathcal{U}_{A^{t}}) & \displaystyle ~\cong ~\bigoplus_{\substack{ n \text{~ even~} \\ 0 \leq n \leq d}}coker(1-B_{n}), \textrm{~and~} \\
K_{1}(\mathcal{U}_{A^{t}}) & ~\cong~ \bigoplus_{\substack{n \text{~odd~} \\ 0 \leq n \leq d}}coker(1-B_{n}) \oplus \mathbb{Z}.
\end{align*}
\item If $d$ is even and $\det(A)<0$ then
\begin{align*}
 K_{0}(\mathcal{U}_{A^{t}}) & \displaystyle ~\cong ~\bigoplus_{\substack{ n \text{~ even~} \\ 0 \leq n \leq d}}coker(1-B_{n}), \textrm{~and~} \\
K_{1}(\mathcal{U}_{A^{t}}) & ~\cong~ \bigoplus_{\substack{n \text{~odd~} \\ 0 \leq n \leq d}}coker(1-B_{n}). 
\end{align*}
\item If $d$ is odd and $\det(A)>0$ then 
\begin{align*}
 K_{0}(\mathcal{U}_{A^{t}}) & \displaystyle ~\cong ~\bigoplus_{\substack{ n \text{~ even~} \\ 0 \leq n \leq d}}coker(1-B_{n}) \oplus \mathbb{Z}, \textrm{~and~}\\
K_{1}(\mathcal{U}_{A^{t}}) & ~\cong~ \bigoplus_{\substack{n \text{~odd~} \\ 0 \leq n \leq d}}coker(1-B_{n}). 
\end{align*}
\item If $d$ is odd and $\det(A)<0$ then
\begin{align*}
 K_{0}(\mathcal{U}_{A^{t}}) & \displaystyle ~\cong ~\bigoplus_{\substack{ n \text{~ even~} \\ 0 \leq n \leq d}}coker(1-B_{n}), \textrm{~and~} \\
K_{1}(\mathcal{U}_{A^{t}}) & ~\cong~ \bigoplus_{\substack{n \text{~odd~} \\ 0 \leq n \leq d}}coker(1-B_{n}). 
\end{align*}
\end{enumerate}
\end{thm}

\bibliography{references.bib}

\def\cprime{$'$} \def\cprime{$'$} \def\cprime{$'$}
\providecommand{\bysame}{\leavevmode\hbox to3em{\hrulefill}\thinspace}
\providecommand{\MR}{\relax\ifhmode\unskip\space\fi MR }
\providecommand{\MRhref}[2]{%
  \href{http://www.ams.org/mathscinet-getitem?mr=#1}{#2}
}
\providecommand{\href}[2]{#2}
\begin{thebibliography}{EaHR10}

\bibitem[BE10]{Ex2}
Giuliano Boavo and Ruy Exel, \emph{Partial crossed product description of the
  {$C^*$}-algebras associated to integral domains}, arxiv:1010.0967v2/math.OA,
  2010.

\bibitem[Bla87]{Bla}
B.~Blackadar, \emph{K-theory for operator algebras}, Springer Verlag,Newyork,
  1987.

\bibitem[CL09]{Cuntz-Li2}
Joachim Cuntz and Xin Li, \emph{K-theory of ring {$C^{*}$}-algebras associated
  to function fields}, arxiv:0911.5023v1, 2009.

\bibitem[CL10]{Cuntz-Li}
Joachim Cuntz and Xin Li, \emph{The regular {$C^\ast$}-algebra of an integral
  domain}, Quanta of maths, Clay Math. Proc., vol.~11, Amer. Math. Soc.,
  Providence, RI, 2010, pp.~149--170. \MR{2732050}

\bibitem[CL11]{Cuntz-Li-1}
\bysame, \emph{{$C^\ast$}-algebras associated with integral domains and crossed
  products by actions on adele spaces}, J. Noncommut. Geom. \textbf{5} (2011),
  no.~1, 1--37. \MR{2746649}

\bibitem[Cun08]{Cuntz}
Joachim Cuntz, \emph{{$C^*$}-algebras associated with the {$ax+b$}-semigroup
  over {$\Bbb N$}}, {$K$}-theory and noncommutative geometry, EMS Ser. Congr.
  Rep., Eur. Math. Soc., Z\"urich, 2008, pp.~201--215. \MR{2513338
  (2010i:46086)}

\bibitem[EaHR10]{Exel-Rae}
Ruy Exel, Astrid an~Huef, and Iain Raeburn, \emph{Purely infinite simple
  ${C}^{*}$-algebras associated to integer dilation matrices},
  arxiv.1003.2097/math.OA, 2010.

\bibitem[FS81]{Fack-Skandalis}
Thierry Fack and Georges Skandalis, \emph{Connes' analogue of the {T}hom
  isomorphism for the {K}asparov groups}, Invent. Math. \textbf{64} (1981),
  no.~1, 7--14. \MR{621767 (82g:46113)}

\bibitem[KLQ11]{Quigg-Landstad}
S.~Kaliszewski, M.~Landstad, and J.~Quigg, \emph{A crossed-product approach to
  the {C}untz-{L}i algebras}, arxiv:1012:5285v3, 2011.

\bibitem[LL12]{wolf-Li}
Wolfgang Luck and Xin Li, \emph{K-theory for ring ${C}^{*}$ algebras- the case
  of number fields with higher roots of unity}, arxiv/1201.4296, 2012.

\bibitem[MRR]{Rae-Laca-Rammage}
Laca Marcelo, Iain Raeburn, and Jacqui Rammage, \emph{Phase transition of
  {E}xel crossed products associated to dilation matrices},
  arxiv/math.OA:1101.4713v1.

\bibitem[Sun12]{Sundar-Cuntz-Li}
S.~Sundar, \emph{Cuntz-{L}i relations, {I}nverse semigroups and {G}roupoids},
  arxiv:1201.4620v1, 2012.

\end{thebibliography}
\bibliographystyle{amsalpha}

\nocite{Bla}

\end{document}